# Neighbor-Gradient Single-Pass Method for solving anisotropic eikonal equation

Myong-Song Ho[1*], Ji-Song Pak[1], Song-Mi Jo[1], Ju-Song Kim[1]

[1] Faculty of Mathematics, Kim Il Sung University, Pyongyang, DPR Korea

[*] Corresponding Author: Myong-Song Ho

Email: ms.ho0511@ryongnamsan.edu.kp

Conflict of Interest: The authors declare that they have no conflict of interest.

Funding: This research did not receive any specific grant from funding agencies in the public, commercial, or not-for-profit sectors.

**Abstract**

We develop a single pass method for approximating the solution to an anisotropic eikonal equation related to the anisotropic min-time optimal trajectory problem. Ordered Upwind Method (OUM) solves this equation, which is a single-pass method with an asymptotic complexity. OUM uses the search along the accepted front (SAAF) to update the value at considered nodes. Our technique, which we refer to as "Neighbor-Gradient Single-Pass Method", uses the minimizer of the Hamiltonian, in which the gradient is substituted with neighbor gradient information, to avoid SAAF. Our technique is considered in the context of control-theoretic problem. We begin by discussing SAAF of OUM. We then prove some properties of the value function and its gradient, which provide the key motivation for constructing our method. Based on these discussions, we present a new single-pass method, which is fast since it does not require SAAF. We test this method with several anisotropic eikonal equations to observe that it works well while significantly reducing the computational cost.

## 1. Introduction

Hamilton-Jacobi equations are considered in several contexts, including classical mechanics, front propagation, control problems and differential games. For optimal control problems, In particular, the value function can be characterized as the unique viscosity solution of a Hamilton–Jacobi–Bellman equation.

In this paper, we deal with generalized eikonal equations, which are associated to some min-time optimal problems. These equations are used not only in optimal control problems but also in many fields, such as seismic travel time estimation, geodesics problems, etc.

Some types of generalized eikonal equations associated to some min-time optimal problems are given as follows ([1]):

$$\sup_{a \in S_1}\{-a \cdot \nabla u(x)\} = 1 \quad \text{(homogeneous eikonal equation)} \tag{1.1}$$

$$\sup_{a \in S_1}\{-f_1(x) a \cdot \nabla u(x)\} = 1 \quad \text{(nonhomogeneous eikonal equation)} \tag{1.2}$$

$$\sup_{a \in S_1}\{-f_2(a) a \cdot \nabla u(x)\} = 1 \quad \text{(homogeneous anisotropic eikonal equation)} \tag{1.3}$$

$$\sup_{a \in S_1}\{-f_3(x, a) a \cdot \nabla u(x)\} = 1 \quad \text{(nonhomogeneous anisotropic eikonal equation)} \tag{1.4}$$

where $x \in \Omega \setminus \partial\Omega$, $\Omega$ is a bounded convex open set in $R^d$, $\partial\Omega$ is a nonempty target set and $f_1$, $f_2$, $f_3$ are given positive and Lipchitz continuous scalar functions, and $S_1$ is the unit ball in $\mathbf{R}^d$, representing the set of the admissible controls. We will always consider the Dirichlet condition $u = 0$ on $\partial\Omega$. To simplify the notations, we restrict the discussion to the case $d = 2$.

Various numerical methods based on Semi-Lagrangian schemes (which directly discretize Bellman's dynamic programming principle) for solving equations (1.1) -(1.4) have been proposed. ([2-5])

Fast marching method (FMM) ([6-9]) has been proposed to solve the isotropic eikonal equations such as (1.1) and (1.2). FMM is a local single-pass method, which is extremely fast.

An algorithm is said to be *single-pass* if each grid node is recomputed at γ times where γ is a priori a (small) number which depends only on the equation and on the grid structure (not on the number of grid nodes) and said to be *local* if the computation at any grid node involves only the values of first neighboring nodes([[1, 10]). Single-pass methods divide the grid nodes in three groups: Accepted ($Acc$) region, Considered ($Cons$) region and *Far* (Far) region. $Acc$ is the set of grid nodes whose values are definitively computed, $Cons$ is the set of grid nodes whose values are computed but their values are not yet final, and *Far* is the set of grid nodes whose values are not yet computed. ([1, 10])

Unfortunately, FMM fails for anisotropic eikonal equations such as (1.3) and (1.4) because the neighboring simplex from which the characteristic approaches a node $x$ may contain another node $y$ such that causality does not hold: $u(x) \leq u(y)$.

For the anisotropic eikonal equation, Ordered Upwind Method (OUM), which is a single-pass method, has been proposed ([11, 12]). OUM can solve the general anisotropic eikonal equation and it has an asymptotic complexity only a constant factor (related to the degree of anisotropy) worse than FMM. It is shown that the solutions produced by OUM in the boundary value formulation converge at a rate of at least the square root of the largest edge length in the mesh in terms of maximum error. ([13]). OUM avoids the difficulty of FMM by searching along the accepted front (SAAF) to find a set of neighboring nodes (which may not be direct neighbors) whose values have been accepted, and then constructing a virtual simplex with these nodes from which to update. Although SAAF does not degrade the asymptotic complexity of OUM, it does significantly increase the computational cost in practice ([14-16]). MAOUM (Monotone Acceptance OUM), which is a generalization of OUM, has been proposed. MAOUM computes the solution in a Fast Marching fashion, but employs large stencils that are pre-computed for each grid, and thus MAOUM is a two-pass method. ([1])

Efficient methods such as Fast Sweeping Method (FSM) ([17-22]) and Fast Iterative Method (FIM) ([23-25]) have been proposed for solving anisotropic eikonal equations. These are not single-pass method.

On the other hand, for solving the generalized eikonal equation in a moving medium, Characteristic Fast Marching Method (CFMM) ([26, 27]), Operator Split-based Method ([28, 29]) and Improved Characteristic Fast Marching Method (ICFMM) ([30]) have been proposed.

In this paper, we focus on OUM which is a single-pass method for solving the anisotropic eikonal equation (1.3) and (1.4).

We need to recall the procedures for updating the value at $x \in Cons$ in OUM:

(a) Search along the accepted front to find the neighboring nodes (which may not be direct neighbors of $x$) that the distance from the line segment joining them to $x$ is smaller than $h\gamma$, where $h$ is diameter of grid and $\gamma$ is the anisotropic coefficient defined in [11].

(b) For every found neighboring nodes, construct a virtual simplex from these nodes and calculate the approximated value from this simplex using the semi-Lagrangian scheme.

(c) Update the value at $x$ with the minimum value among the approximated values for every found neighboring nodes.

Here, procedure (a) is SAAF which increases the computational cost. There exists one simplex which contains the characteristic for $x$ among the virtual simplexes constructed by procedure (b). So procedures (a), (b) and (c) allow to update the value at $x$ from the characteristic direction.

If the characteristic direction at $x$ was a priori known, the computational cost decreases, since one would be able to update the value from the characteristic direction, without SAAF .i.e., update the value from the intersection point between the accepted front and the half line( initiated at the point $x$ and pointing in the characteristic direction). Since the anisotropic eikonal equations (1.3) and (1.4) are nonlinear, one cannot know the characteristic direction in advance.

We also need to recall the following well-known formula ([11]):

$$a_1 = \arg\min_{a \in S_1}\{(\nabla u(x) \cdot a) f(x,a)\} \quad (1.5)$$

, where $a_1$ is the optimal control value (the characteristic direction) at $x$ (i.e., the minimizer of Hamiltonian). From this formula, if $\nabla u(x)$ was a priori known, one would be able to obtain the optimal control value according to (1.5) and update the value at $x$ without using SAAF. Unfortunately, $\nabla u(x)$ is also unknown unless $x$ is accepted.

We use a technique for updating the value at $x$ from the minimizer of Hamiltonian, in which the gradient is substituted with the gradient at a neighboring accepted node.

The following questions may arise in using the above technique.

First question is that whether there exists at least one accepted neighboring grid node such that the gradient at this node is similar to $\nabla u(x)$.

Second question is that whether the resulting error caused by using the minimizer of Hamiltonian (substituted with the gradient at the above neighboring accepted node) instead of the optimal control value may accumulate during the iteration of the algorithm and ultimately may fail to ensure convergence.

The answer to the first question is given in Lemma 2.1 of Section 2.

Before answering the second question, we pay attention to the following fact.

In the case of using OUM, one assumes that the approximated values at the accepted nodes are the same as the real values of value function, at a certain stage of the algorithm. Suppose that $\tilde{x}$ is the node with smallest value among all the considered nodes, then the error between the approximated value and real value at $\tilde{x}$ is $O(h^2)$. This fact was demonstrated for validating OUM in [11]. The answer to the second question is similar as suggested above. Detailed proof is given in Sections 2 and 4. We call the resulting method the *Neighbor-Gradient Single-Pass Method*.

In Section 2, we prove some properties related to the validity of our technique, including the answers to two above questions. Based on this, in Section 3, we present approximation scheme that updates the values at considered nodes by using the gradients at accepted neighboring nodes. In Section 4, we present a neighbor-gradient single-pass algorithm, implementation and computational complexity of our algorithm. In Section 5, we test the efficiency and convergence through numerical experiments. Section 6 discusses future directions and concludes the paper. Some basic properties of the value function are also presented in Appendix.

## 2. Anisotropic min-time problem and Properties of the value function.

### 2.1 Statement of problem

We consider an optimal trajectory problem for a vehicle moving with the speed depending upon the direction of motion and the current position of the vehicle inside the domain $\Omega$. ([11])

The dynamics of the vehicle is defined by

$$y'(t) = f(y(t), \alpha(t))a(t)$$
$$y(0) = x \in \Omega$$

, where $y(t)$ is the position of the vehicle at time $t$, $S_1 = \{a \in \mathbf{R}^2 \mid \|a\| = 1\}$ is the set of admissible control values, and $A = \{\alpha : \mathbf{R}_{+,0} \mapsto S_1 \mid \alpha(\cdot) \text{ is measurable}\}$ is the set of admissible controls.

We are interested in studying $y(t)$ only while the vehicle remains inside $\Omega$, i.e., until the exit time

$$T(x, \alpha(\cdot)) = \inf\{t \in \mathbf{R}_{+,0} \mid y(t) \in \partial\Omega\}.$$

Then a min-time optimal trajectory problem is to find a control $\alpha(\cdot)$ which minimizes the exit time $T(x, \alpha(\cdot))$.

The value function of the min-time optimal trajectory problem is defined by

$$u(x) = \inf_{\alpha(\cdot) \in A} T(x, \alpha), \quad x \in \Omega \quad (2.1)$$

Then this value function is the unique viscosity solution of the following Hamilton-Jacobi-Bellman's equation. ([11])

$$\min_{a \in S_1}\{f(x, a) a \cdot \nabla u(x)\} = -1, \quad x \in \Omega \setminus \partial\Omega$$
$$u(x) = 0, \quad x \in \partial\Omega \quad (2.2)$$

It is clear that equation (2.1) coincides the anisotropic eikonal equation (1.4).

We will assume that there exists constants $F_1, F_2$ such that

$$F_1 \leq f(x, a) \leq F_2 \quad (F_1, F_2 > 0) \text{ for all } x \in \Omega \text{ and } a \in S_1.$$

and that $f(x, a)$ is Lipschitz continuous in the state variable, uniformly in the control value, that is

$$|f(x, a) - f(y, a)| < L\|x - y\|, \text{ for all } x, y \in \Omega, \ a \in A.$$

Hamilton-Jacobi-Bellman PDE (2.1) can be rewritten in the form $H(x, \nabla u(x)) = -1$, where the Hamiltonian $H(x, p) = -\min_{a \in S_1}\{p \cdot a f(x, a)\}$. Moreover $H$ is convex in the second argument. ([11])

Assume that a regular rectangular grid $X$ is defined on $\Omega$, where the grid sizes in $x$- and $y$-coordinates are the same .i.e., $\Delta x = \Delta y$.

We define $h := \sqrt{2}\Delta x$ and for each grid node $x = (x_i, y_j) \in X$, we denote

$$N(x) := \{(x_{i-1}, y_j), (x_{i+1}, y_j), (x_i, y_{j-1}), (x_i, y_{j+1})\}$$

$$ND(x) := \{(x_{i-1}, y_j), (x_{i+1}, y_j), (x_i, y_{j-1}), (x_i, y_{j+1}), (x_{i-1}, y_{j-1}), (x_{i-1}, y_{j+1}), (x_{i+1}, y_{j-1}), (x_{i+1}, y_{j+1})\}.$$

## 2.2 The properties of the value function

We will prove some properties of the value function (2.2) and its gradient.

**Lemma 2.1.** *Consider two points* $x_1, x_2 \in \Omega$ *and assume that* $u$ *is differentiable at* $x_1$. *We define*

$$V(x_1, x_2) := u(x_2) + \frac{\|x_2 - x_1\|}{f(x_1, a_1)} \qquad (2.3)$$

, *where* $a_1 = x_2 - x_1 / \|x_2 - x_1\|$. *Then*

$$V(x_1, x_2) - u(x_1) = \|x_1 - x_2\| \left( \nabla u(x_1) \cdot a_1 + \frac{1}{f(x_1, a_1)} \right) + O(\|x_1 - x_2\|^2).$$

*Proof.* By the definition of $V(x_1, x_2)$ and differentiability of $u$ at $x_1$, we obtain

$$V(x_1, x_2) = u(x_1) + u(x_2) - u(x_1) + \frac{\|x_1 - x_2\|}{f(x_1, a_1)} = u(x_1) + \nabla u(x_1) \cdot a_1 \|x_1 - x_2\| +$$

$$+ O(\|x_1 - x_2\|^2) + \frac{\|x_1 - x_2\|}{f(x_1, a_1)} = u(x_1) + \|x_1 - x_2\| \left( \nabla u(x_1) \cdot a_1 + \frac{1}{f(x_1, a_1)} \right) + O(\|x_1 - x_2\|^2). \square$$

**Lemma 2.2.** *Consider a grid node* $x \in \Omega$. *Suppose that the characteristic for* $x$ *intersects the line segment* $x_1 x_2$ *at* $x_0$, *where* $x_1, x_2 \in ND(x)$ *and* $x_1 \in N(x_2)$. *Assume that* $u$ *is differentiable at* $x$, $x_1$ *and* $x_2$. *Then*

$$\|\nabla u(x_i) - \nabla u(x)\| = O(h), i = 1, 2.$$

*Proof.* By Lemma 3 in Appendix, we get

$$\|\nabla u(x_0) - \nabla u(x)\| = O(h). \quad (2.4)$$

Without loss of generality, let $x_2 \in N(x)$, $x_1 \in ND(x)(x_1 \notin N(x))$.

(i) Let us prove that $\|\nabla u(x_2) - \nabla u(x)\| = O(h)$.

By the fact that characteristics never emanate from the shocks-non differentiable point ([11]) and that $u$ is differentiable at $x_2$, we see that $u$ is differentiable at $x_0$.

Therefore, we get

$$u(x_2) = u(x_0 + x_2 - x_0) = u(x_0) + u_x(x_0)\|x_2 - x_0\| + O\|x_2 - x_0\|^2$$
$$u(x_0) = u(x_2 + x_0 - x_2) = u(x_2) - u_x(x_2)\|x_2 - x_0\| + O\|x_2 - x_0\|^2$$

and thus

$$(u_x(x_0) - u_x(x_2))\|x_2 - x_0\| = O\|x_2 - x_0\|^2$$
$$u_x(x_0) - u_x(x_2) = O\|x_2 - x_0\| = O(h). \qquad (2.5)$$

Using (2.4) and (2.5), it's clear that

$$u_x(x) - u_x(x_2) = O(h) \qquad . \qquad (2.6)$$

Similarly, we can get

$$u_y(x) - u_y(x_2) = O\|x_2 - x\| = O(h) . \qquad (2.7)$$

Therefore, from (2.6) and (2.7), we obtain $\|\nabla u(x) - \nabla u(x_2)\| = O(h)$.

(ii) Let us prove that $\|\nabla u(x_1) - \nabla u(x)\| = O(h)$.

Similar to the proof of (i), we can get $u_x(x_1) - u_x(x_0) = O(h)$.

Using (2.4), we obtain $u_x(x_0) - u_x(x) = O(h)$. Therefore,

$$u_x(x_1) - u_x(x) = O(h) . \qquad (2.8)$$

Considering the directional derivative of $u$ at $x_1$ in the direction $\overrightarrow{x_1 x}$, we can get

$$u_{\overrightarrow{x_1 x}}(x_1) = \frac{1}{\sqrt{2}} u_x(x_1) + \frac{1}{\sqrt{2}} u_y(x_1) = \frac{u(x) - u(x_1)}{\|x - x_1\|} = \frac{u(x_0) + \|x - x_0\|/f(x_0, a) - u(x_1)}{\sqrt{2}h}$$

and thus

$$u_x(x_1) + u_y(x_1) = \frac{u(x_0) - u(x_1)}{h} + \frac{\|x - x_0\|}{f(x_0, a)h} = \frac{u(x_0) - u(x_1)}{\|x_0 - x_1\|} \frac{\|x_0 - x_1\|}{h} + \frac{\|x - x_0\|}{f(x_0, a)h} =$$

$$= u_x(x_1) \frac{\|x_0 - x_1\|}{h} + \frac{\|x - x_0\|}{f(x_0, a)h} + O(h).$$

Therefore,

$$u_x(x_1) \frac{\|x_0 - x_2\|}{h} + u_y(x_1) = \frac{\|x - x_0\|}{f(x_0, a)h} + O(h) \qquad . \qquad (2.9)$$

Similarly, we can get

$$u_x(x) \frac{\|x_0 - x_2\|}{\|x - x_0\|} + u_y(x) \frac{h}{\|x - x_0\|} = \frac{u(x) - u(x_0)}{\|x - x_0\|} = \frac{1}{f(x_0, a)}$$

$$u_x(x) \frac{\|x_0 - x_2\|}{h} + u_y(x) = \frac{\|x - x_0\|}{f(x_0, a)h} \qquad . \qquad (2.10)$$

Then from (2.9) and (2.10), it's clear that

$$u_y(x_1) - u_y(x) = (u_x(x) - u_x(x_1)) \frac{\|x_0 - x_2\|}{h} = \frac{\|x_0 - x_2\|}{h} O(h) = O(h). \qquad (2.11)$$

Using (2.8) and (2.11), we obtain $\|\nabla u(x) - \nabla u(x_1)\| = O(h)$. □

**Lemma 2.3** *Consider a grid node $x \in \Omega$. Suppose that the characteristic for $x$ intersects the line segment $x_1x_2$ at $x_0$, where $x_1, x_2 \in ND(x)$ and $x_1 \in N(x_2)$. Assume that $u$ is differentiable at $x, x_1$ and $x_2$. Suppose that $a_0$ is the optimal control value at $x$. For $x' \in \{x_1, x_2\}$, we denote $\hat{a} := \arg\min_{a \in S_1}\{\nabla u(x')af(x,a)\}$ Then*

$$\|\hat{a} - a_0\| = O(h)$$

$$\nabla u(x) \cdot \hat{a} + \frac{1}{f(x,\hat{a})} = O(h).$$

*Proof.* Let us prove that $\nabla u(x) \cdot \hat{a} + \frac{1}{f(x,\hat{a})} = O(h)$.

Since $a_0$ is the optimal control value at $x$,

$$a_0 = \arg\min_{a \in S_1}\{\nabla u(x)af(x,a)\}, \quad \nabla u(x)a_0 f(x,a_0) = -1. \quad (2.12)$$

By Lemma 2.2, we get $\|\nabla u(x) - \nabla u(x')\| = O(h)$ and thus

$$|\nabla u(x) \cdot a_0 f(x,a_0) - \nabla u(x') \cdot a_0 f(x,a_0)| = |(\nabla u(x) - \nabla u(x')) \cdot a_0 f(x,a_0)|$$
$$\leq \|\nabla u(x) - \nabla u(x')\| f(x,a_0) = O(h).$$

Therefore,

$$\nabla u(x') \cdot a_0 f(x,a_0) = \nabla u(x) \cdot a_0 f(x,a_0) + O(h) = -1 + O(h).$$

Since $\hat{a} = \arg\min_{a \in S_1}\{\nabla u(x')af(x,a)\}$, the following inequality holds:

$$\nabla u(x') \cdot \hat{a} f(x,\hat{a}) < \nabla u(x') \cdot a_0 f(x,a_0) = -1 + O(h).$$

By Lemma 2.2, $\|\nabla u(x) - \nabla u(x')\| = O(h)$. Therefore,

$$\nabla u(x) \cdot \hat{a} f(x,\hat{a}) = \nabla u(x') \cdot \hat{a} f(x,\hat{a}) + O(h) < -1 + O(h).$$

Using (2.12), we obtain $\nabla u(x) \hat{a} f(x,\hat{a}) \geq -1$.

Since $\nabla u(x) \cdot \hat{a} f(x,\hat{a}) = -1 + O(h)$ and $f$ is upper bound, we get

$$\nabla u(x) \cdot \hat{a} + \frac{1}{f(x,\hat{a})} = O(h).$$

Let us prove that $\|\hat{a} - a_0\| = O(h)$.

We define $A(p) := \arg\min_{a \in S_1}\{paf(x,a)\}$. In particular, we set $A(\nabla u(x')) = \hat{a}$.

We denote $H'(p) := \nabla u(x) A(p) f(x, A(p))$. Then $H'(p)$ is continuous.

Let us consider the sequence $\{p_n\}$ such that $p_n \to p_0 (:= \nabla u(x))$.

Since $A(p_n)$ is bounded sequence, it has convergent subsequence $\{A(p_k)\}$ ($A(p_k) \to a^*$).

By the fact that $p_k \to p_0$ and that $H'(p)$ is continuous,

$$H'(p_k) = \nabla u(x) A(p_k) f(x, A(p_k)) \to H'(\nabla u(x)) = -1.$$

On the other hand, by the fact that $\nabla u(x) a f(x,a)$ is continuous in $a$, we get

$$\nabla u(x) A(p_k) f(x, A(p_k)) \to \nabla u(x) a^* f(x, a^*)$$

and $\nabla u(x) a^* f(x, a^*) = -1$. Therefore $a^* = a_0$.

Since $A(p_n)$ is bounded sequence and the limit of all convergent subsequences are the same as $a_0$, we get $A(p_n) \to a_0$.

Therefore, $A(p)$ is continuous at $p_0 = \nabla u(x)$.

To complete the proof, we recall that $A(\nabla u(x')) = \hat{a}$, $A(\nabla u(x)) = a$ and $\|\nabla u(x) - \nabla u(x')\| = O(h)$. □

Form the above lemmas (lemma 2.1, 2.2 and 2.3), we can claim the following theorem which provides the key motivation for constructing our method.

**Theorem 2.1** *Consider a simple closed curve* $\Gamma \subset \Omega \setminus \partial \Omega$ *with the property that for any point* $x$ *on* $\Gamma$ *there exists a gird node* $\hat{x}$ *inside* $\Gamma$ *such that* $\|x - \hat{x}\| < h$. *Suppose a grid node* $\bar{x}$ *is such that* $u(\bar{x}) \le u(x_i)$ *for all the grid nodes* $x_i \in X$ *inside* $\Gamma$. *Then the following statements hold.*

*(i) For the grid node* $\bar{x}$, *Suppose that the characteristic for* $\bar{x}$ *intersects the line segment* $x_1 x_2$ *at* $x_0$, *where* $x_1, x_2 \in ND(\bar{x})$ *and* $x_1 \in N(x_2)$. *Then at least one of* $x_1$ *and* $x_2$ *is on* $\Gamma$ *or outside* $\Gamma$.

*(ii) Suppose that* $x_1$ *is the grid node on* $\Gamma$ *or outside* $\Gamma$ *and we denote*

$$\hat{a} := \arg\min_{a \in S_1} \{\nabla u(x_1) a f(\bar{x}, a)\}, \quad V(\bar{x}, \hat{a}) := u(\hat{x}) + \frac{\|\bar{x} - \hat{x}\|}{f(\bar{x}, \hat{a})}.$$

*Then*

$$V(\bar{x}, \hat{a}) - u(\bar{x}) = O(h^2),$$

*where* $\hat{x}$ *is the intersection point between* $\Gamma$ *and a half-line (initiated at the point* $\bar{x}$ *and pointing in the direction of* $\hat{a}$.).

*Proof.* (i) By Lemma 4 in Appendix, we get that $\min(u(x_1), u(x_2)) < u(\bar{x})$.

Without loss of generality, assume that $u(x_1) < u(\bar{x})$.

Since $u(\bar{x}) \leq u(x_i)$ for all the grid nodes $x_i \in X$ inside $\Gamma$, $x_1$ is not inside $\Gamma$. Therefore $x_1$ is on $\Gamma$ or outside $\Gamma$.

(ii) Suppose that $a_0$ is the optimal control value at $\bar{x}$. i.e., $a_0 = \arg\min_{a \in S_1}\{\nabla u(\bar{x}) a f(\bar{x}, a)\}$.

By [11, Lemma 3.1], we get

$$\|\hat{x} - \bar{x}\| \leq h\gamma \qquad (2.13)$$

We define $\hat{a} := \dfrac{\hat{x} - \bar{x}}{\|\hat{x} - \bar{x}\|}$. By the Mean Value Theorem of Differentially, we get $\hat{a} \wedge a_0 = O(h)$.

On the other hand, By Lemma 2.3, we obtain $\hat{a} \wedge a_0 = O(h)$. Therefore, $\hat{a} \wedge \hat{a} = O(h)$.

Using the continuity of $\Gamma$, we get $\|\tilde{x} - \bar{x}\| - \|\hat{x} - \bar{x}\| = O(h)$.

Combining with (2.13), we get $\|\hat{x} - \bar{x}\| = O(h)$.

By Lemma 2.1, we obtain

$$V(\bar{x}, \hat{a}) - u(\bar{x}) = \|\hat{x} - \bar{x}\|\left(\nabla u(\bar{x}) \cdot \hat{a} + \frac{1}{f(\bar{x}, \hat{a})}\right) + O(\|\hat{x} - \bar{x}\|)^2.$$

By the fact that $\|\hat{x} - \bar{x}\| = O(h)$ and Lemma 2.3, we get

$$\nabla u(\bar{x}) \cdot \hat{a} + \frac{1}{f(\bar{x}, \hat{a})} = O(h).$$

Therefore, $V(\bar{x}, \hat{a}) - u(\bar{x}) = O(h^2)$. □

*Remark 2.1* When we consider this theorem in the context of algorithm, $\Gamma$ is accepted front and the fact that $x'$ is on $\Gamma$ or outside $\Gamma$ means that $x'$ is in $Acc$, therefore it is possible to approximate the gradient of value function at $x'$.

### 3. Approximation of the value function

As before, all the grid nodes are divided into 3 groups: $Acc$ (accepted region), $Cons$ (considered region), $Far$ (far region). $AFF$ is defined as a set of the accepted grid nodes $x_i$ such that there exists a considered node $\tilde{x}(\in ND(x_i))$. Define the set $AF$ of the line segments $x_j x_k$, where $x_j, x_k$ are grid nodes in the $AFF$ such that $x_j \in N(x_k)$. For considered grid node $x$, we denote

$$NF(x) := \{x_j x_k \in S_{AF} \mid \exists \tilde{x} \in x_j x_k, \|\tilde{x} - x\| < h\gamma\}.$$

Based on Theorem 2.1, we use a technique for approximating the value function by using a gradient at accepted neighboring grid nodes.

We define $U$ to be the definitely computed value at accepted grid nodes.

For the grid node $\hat{x} = (x_i, y_j) \in Acc$, we obtain $\nabla U(\hat{x}) = (DX, DY)$ as follows.

Let $x_L = (x_{i-1}, y_j)$, $x_R = (x_{i+1}, y_j)$, $x_D = (x_i, y_{j-1})$, $x_U = (x_i, y_{j+1})$.

if $x_L \in Acc$, $x_R \in Acc$ then $DX = (U(x_R) - U(x_L))/(2\Delta x)$

if $x_L \in Acc$, $x_R \notin Acc$ then $DX = (U(\hat{x}) - U(x_L))/\Delta x$

if $x_R \in Acc$, $x_L \notin Acc$ then $DX = (U(x_R) - U(\hat{x}))/\Delta x$

if $x_R \notin Acc$, $x_L \notin Acc$ then $DX = 0$

if $x_D \in Acc$, $x_U \in Acc$ then $DY = (U(x_U) - U(x_D))/(2\Delta y)$

if $x_D \in Acc$, $x_U \notin Acc$ then $DY = (U(\hat{x}) - U(x_D))/\Delta y$

if $x_U \in Acc$, $x_D \notin Acc$ then $DY = (U(x_U) - U(\hat{x}))/\Delta y$

if $x_U \notin Acc$, $x_D \notin Acc$ then $DY = 0$

Now we propose a way to obtain the approximate value of $U$ at $x \in Cons$ using the gradients at the neighboring grid nodes belonging to $Acc$. For $x \in Cons$, we define

$$G(x) := \bigcup_{\bar{x} \in ND(x) \cap Acc} \nabla U(\bar{x}), \quad A(x) := \left\{ \hat{a} = \arg\min_{a \in S_1} \{\tilde{g} \cdot af(x,a), \ \tilde{g} \in G(x)\right\}.$$

Let $l_{\hat{a}}$ be a half line initiated at the point $\mathbf{x}$ and pointing in the direction of $\hat{a}$.

Suppose that the half line $l_{\hat{a}}$ intersects the segment $x_{j_0} x_{k_0} \in NF(x)$ at the point $\tilde{x}$. We set

$$V_{\hat{a}}(x) := V_{x_{j_0}, x_{k_0}, \tilde{x}}(x) = \frac{\|\tilde{x} - x\|}{f(x, \tilde{x} - x/\|\tilde{x} - x\|)} + \lambda U(x_{j_0}) + (1-\lambda)U(x_{k_0}),$$

where $\tilde{x} = \lambda x_{j_0} + (1-\lambda)x_{k_0}$ and if none of every half line $l_{\hat{a}}$ intersect $NF(x)$, put $V_{\hat{a}}(x) = \infty$.

For $x' \in ND(x) \cap Acc$ we define $V_{x'}(x) = \frac{\|x' - x\|}{f(x, x'-x/\|x'-x\|)} + U(x')$.

From the above discussion, we update the value at $x$ as following.

$$V(x) := \{\min_{\hat{a} \in A(x)} V_{\hat{a}}(x), \min_{x' \in ND(x) \cap Acc} V_{x'}(x)\}. \quad (3.1)$$

Here, we introduce the approximation scheme of OUM ([11]).

$$V(x) := \min_{x_j x_k \in NF(x)} \min_{x' \in x_j x_k} V_{x_j x_k x'}(x). \quad (3.2)$$

## 4. Neighbor-Gradient Single Pass Algorithm

### 4.1 Algorithm

1. Move all grid nodes into *Far* and put $U(x)=\infty$ for all $x \in X$.

2. Move all grid nodes $y (\in \partial\Omega \cap X)$ into *Acc*. ($U(y)=0$)

3. Update the value at the *Far* grid nodes neighboring to *Acc* by using (3.2) and move them into *Acc*, and update the gradient at all grid nodes in *Acc*. (See section 3)

4. Move all the *Far* grid nodes neighboring to *Acc* into *Cons* and evaluate the tentative values $V$ at these points by using (3.1).

5. Find the grid node $\bar{x}$ with the smallest value of $V$ among all the grid nodes in *Cons*.

6. Move $\bar{x}$ into *Acc* ($U(\bar{x})=V(\bar{x})$), and reevaluate the gradient at $\bar{x}$ and at all grid nodes in $ND(\bar{x}) \cap Acc$. (See section 3)

7. Reevaluate the tentative value $V$ for all the grid nodes in $ND(\bar{x}) \cap Acc^c$ and move the grid nodes in $ND(\bar{x}) \cap Far$ into *Cons*.

8. If *Cons* is not empty, then go to 5.

**Motivation**: Assume that $U(x)=u(x)$ for all $x \in AF$. Suppose that $\hat{x} \in Cons$ is such that $u(\hat{x}) = \min_{x \in Cons} V(x)$. Then by the algorithm, $U(\hat{x})=V(\hat{x})$. Let $\bar{x} \in Cons$ is such that $u(\bar{x}) = \min_{x \in Cons} u(x)$. Suppose that $x_1, x_2 \in ND(\bar{x})$, $x_1 \in N(x_2)$ is such that the characteristic for $\bar{x}$ intersects the line segment $x_1 x_2$. Then, by Theorem 2.1, at least one grid node $\tilde{x}$ of $x_1$ and $x_2$ is in *Acc*. We define $\tilde{a} := \arg\min_{a \in S_1}\{\nabla u(\tilde{x}) a f(x,a)\}$ and suppose that $\hat{x}$ is the intersection point between a half line (initiated at the point $\bar{x}$ and pointing in the direction of $\tilde{a}$) and $AF$. If $\|\hat{x}-\bar{x}\| < \gamma h$, then by Theorem 2.1 and our algorithm, we get $V(\bar{x}) < u(\bar{x}) + O(h^2)$. By the fact that $\hat{x} = \arg\min_{x \in Cons} V(x)$, we obtain $u(\hat{x}) + O(h^2) \leq U(\hat{x}) = V(\hat{x}) \leq V(\bar{x}) \leq u(\bar{x}) + O(h^2)$. Since $u(\bar{x}) \leq u(\hat{x})$, $u(\hat{x}) + O(h^2) \leq U(\hat{x}) \leq u(\bar{x}) + O(h^2) \leq u(\hat{x}) + O(h^2)$. Therefore, $U(\hat{x}) = u(\hat{x}) + O(h^2)$. On the other hand, if $\|\hat{x}-\bar{x}\| > \gamma h$, then By Theorem 2.1, we can see that $u(\hat{x}) + \dfrac{\|\bar{x}-\hat{x}\|}{f(\bar{x},(\hat{x}-\bar{x})/\|\bar{x}-\hat{x}\|)} = u(\bar{x})$ $+ O(h^2)$. Suppose that $\hat{x}$ lies on a line segment $x_1 x_2 (\in AF)$. By the semi-concavity of $u$ (Lemma 1 in Appendix) $\lambda u(x_1) + (1-\lambda)u(x_2) = u(\hat{x}) + O(h^2)$. Without loss of generality, suppose that $u(x_1) \leq u(\hat{x})$ $+ O(h^2)$. By the definition of $AFF$, there exists $\tilde{x} \in Cons$ such that $x_1 \in ND(\tilde{x})$. By our algorithm, we get

$$V(\tilde{x}) \leq u(x_1) + \dfrac{\|x_1-\tilde{x}\|}{f(\tilde{x}, x_1-\tilde{x}/\|x_1-\tilde{x}\|)} \leq u(\bar{x}) + \dfrac{h}{f_1} + O(h^2) \leq u(\bar{x}) + \dfrac{h\gamma}{f_2} + O(h^2) \leq u(\hat{x}) + \dfrac{\|\hat{x}-\bar{x}\|}{f(\bar{x},(\hat{x}-\bar{x})/\|\bar{x}-\hat{x}\|)}$$

$+ O(h^2) = u(\bar{x}) + O(h^2)$ and thus $u(\hat{x}) + O(h^2) \leq U(\hat{x}) = V(\hat{x}) \leq V(\tilde{x}) \leq u(\bar{x}) + O(h^2) \leq u(\hat{x}) + O(h^2)$.

Therefore, $U(\hat{x}) = u(\hat{x}) + O(h^2)$.

*Remark 4.1* In the above argument, we applied Theorem 2.1 under the assumption that $\nabla U(\tilde{x}) = \nabla u(\tilde{x})$, where $\nabla U(\tilde{x}) = (U_x(\tilde{x}), U_y(\tilde{x}))$ is the approximated gradient at $\tilde{x}$ (Section 3) and $\nabla u(\tilde{x}) = (u_x(\tilde{x}), u_y(\tilde{x}))$ is the "real" gradient at $\tilde{x}$. The fact that $\|\nabla u(\tilde{x}) - \nabla u(\bar{x})\| = O(h)$ plays an important role in the proof of Theorem 2.1. Thus, once we prove that $\|\nabla U(\tilde{x}) - \nabla u(\tilde{x})\| = O(h)$, the above argument holds without the assumption that $\nabla U(\tilde{x}) = \nabla u(\tilde{x})$. Let $\tilde{x} = (x_i, y_j)$. If either $(x_{i-1}, y_j)$ or $(x_{i+1}, y_j)$ is in $Acc$, by the gradient approximation scheme in Section 3 and assumption that $U(x) = u(x)$ for all $x \in AF$, $U_x(\tilde{x}) - u_x(\tilde{x}) = O(h)$. On the other hand, If both $(x_{i-1}, y_j)$ and $(x_{i+1}, y_j)$ are not in $Acc$, then by the gradient approximation scheme in Section 3, we get $U_x(\tilde{x}) = 0$. Also, since both $(x_{i-1}, y_j)$ and $(x_{i+1}, y_j)$ are not in $Acc$, we obtain $u(x_{i-1}, y_j) \leq u(\bar{x}), u(x_{i+1}, y) \leq u(\bar{x})$. Therefore $u_x(\tilde{x}) = O(h)$ and $U_x(\tilde{x}) - u_x(\tilde{x}) = O(h)$. Similarly, we can get $U_y(\tilde{x}) - u_y(\tilde{x}) = O(h)$. Therefore $\|\nabla U(\tilde{x}) - \nabla u(\tilde{x})\| = O(h)$.

We present the pseudo code of the algorithm.

In the pseudo code, the function *updategradient*(x) reflects the gradient approximation scheme in Section 3.

```
for each x ∈ X do V(x) ← ∞ and add x to Far
for each x ∈ ∂Ω ∩ X do
  u(x) ← 0
  add x to Acc
  for each neighbor x̄ ∈ Far of x do
      U(x̄) ← min(V(x̄), ‖x̄ − x‖/f (x̄, (x−x̄)/‖x̄−x‖))
      add x̄ to Acc
  end
end
for each x ∈ Acc do UpdateGradient(x)
for each x ∈ Acc do
  for each neighbor x̄ ∈ Far of x do UpdateNarrowBand(x̄)
end
While Cons ≠ ∅ do
  x_min ← argmin_{x∈Cons} V(x)
  add x_min to Acc
  UpdateGradient(x_min)
  for each neighbor x̄_min ∈ Acc of x_min do UpdateGradient(x_min)
  for each neighbor x̄_min ∉ Acc of x_min do UpdateNarrowBand(x_min)
end
```

Table 4.1 Pseudo Code of the Algorithm

```
    if x ∈ Far then remove x from Far and add to Cons
    for each neighbor x̄ ∈ Acc of x do
      U(x) ← min(U(x), update(x, ∇U(x̄)))
      a = argmin_b{∇U(x̄) · bf(x, b)}
      let l be the line segment as {x + at | 0 ≤ t ≤ γΔx}
      let y_1, y_2 ∈ Acc be the first intersect neighboring grid nodes with l
      y' ← l ∩ y_1y_2
      U(x) ← min(U(x), U_{y_1,y_2,y'}(x))
```

Table 4.2 UpdateNarrowBand(x)

## 4.2. Implementation

In the first step of Algorithm, we set $Acc = \partial\Omega \cap X$ and $U(x) = 0 (x \in Acc)$. Therefore one can't approximate the gradient at the nodes in $Acc$. So, we update the value at considered nodes by using OUM and put all the considered nodes into $Acc$. From [11, Lemma 3.1], approximated values at these nodes have second order accuracy. To increase the accuracy of the algorithm, one may combine OUM and our method. In detail, one initially uses OUM until 5% of all grid nodes are accepted. Once 5% of all grid nodes are accepted, one switches to our algorithm with $Acc$ consisting of the accepted nodes.

In step 4 of our algorithm, we use (3.1) in Section 3. There, we have to find the intersection point of $NF(x)$ and $l_a$ without SAAF. Suppose that $l_a$ intersects the line segment $x_{j_0}x_{k_0} (\in NF(x))$. Then, by the definition of $AF$ and $NF(x)$, $x_{j_0}x_{k_0}$ is the line segment with the shortest distance to $x$ among the line segments $x_jx_k (x_j, x_k \in Acc, x_j \in N(x_k))$ intersecting with $l_a$. Let $a = (a_1, a_2)$ and without loss of generality, suppose that $a_1, a_2 > 0$. The half-line $l_a$ intersects with one of the two types of line segments (joining two neighboring nodes). i.e., $(x_k, y_l)(x_k, y_{l+1})$ or $(x_k, y_l)(x_{k+1}, y_l)$. Suppose that $l_a$ intersects with line segment $(x_k, y_l)(x_k, y_{l+1})$ at $(x_k, y'_l)$, where $y_l < y'_l < y_{l+1}$. By the definition of $l_a$, we get $\frac{x_k - x_i}{a_1} = \frac{y'_l - y_i}{a_2}$. Since $x_k - x_i = n\Delta x$, we obtain $x_k = x_i + n\Delta x$, $y'_l = y_j + \frac{a_2}{a_1} n\Delta x$. Therefore $dist(x, (x_k, y'_j)) = \sqrt{a_1^2 + a_2^2} \frac{n\Delta x}{a_1}$. By the definition of $NF(x)$, we get $\sqrt{a_1^2 + a_2^2} \frac{n\Delta x}{a_1} \leq h\gamma$ and thus $n \leq \frac{a_1\gamma\Delta x}{\Delta x\sqrt{a_1^2 + a_2^2}} = \frac{a_1\gamma}{\sqrt{a_1^2 + a_2^2}} \leq \gamma$. Similarly, we can get similar result in the case of that $l_a$ intersects with line segment $(x_k, y_l)(x_{k+1}, y_l)$. From the above results, we can find the line segment $x_{j_0}x_{k_0}$ as follows;

While increasing $n$ from 1 to $\gamma$, stop increasing if $x_n, y_n \in Acc$, where $x_n = \left(x_i + n\Delta x, \left[y_j + \frac{a_2}{a_1} n\Delta x\right]\right)$, $y_n = \left(x_i + n\Delta x, \left[y_j + \frac{a_2}{a_1} n\Delta x\right] + 1\right)$. Suppose that the above increasing stops at $k(1 \leq k \leq \gamma)$. Also, while increasing $n$ from 1 to $\gamma$, stop increasing if $x'_n, y'_n \in Acc$, where $x'_n = \left(\left[x_i + \dfrac{a_1}{a_2}\Delta x\right], y_j + n\Delta x\right)$, $y'_n = \left(\left[x_i + \dfrac{a_1}{a_2}\Delta x\right]+1, y_j + n\Delta x\right)$. Suppose that the above increasing stops at $k'(1 \leq k' \leq \gamma)$. Then the line segment with the shorter distance to $x$ among $x'_{k'} y'_{k'}$ and $x_k y_k$ is just $x_{j_0} x_{k_0}$.

### 4.3. Computational Complexity

Once the grid point $x \in \Omega$ belongs to $Acc$, no computation is needed at the point $x$. For every grid point $x \in \Omega$, the computation at $x$ is performed only when one of its neighboring) grid point is accepted as alive point. A computation at $x$ involves only finding the intersection of a given half line with $NF(x)$. ($O(2\gamma)$ complexity as discussed above, here $\gamma$ is the anisotropy coefficient.). On the other hand, for every grid point, there are at most eight neighboring grid points. Hence, the computational complexity of our method is $O(M \log M + \gamma M)$, where $M$ is the number of grid points and the factor of $\log M$ reflects the necessity to maintain a sorted list of the values of $U$ at the grid points belonging to $Cons$ for step 5 of our algorithm.

### 5. Numerical examples

In this section, we test our method with the following equations ([1, 10, 11]):

| Equation | Dynamics |
|---|---|
| HJB-1 | $f(x, y, a) = m_{\lambda, \mu}(a)a$ |
| HJB-2 | $f(x, y, a) = 1/\sqrt{1 + (\nabla g(x, y) \cdot a)^2}$ |
| HJB-3 | $f(x, y, a) = (1 + |x + y|)m_{\lambda, \mu}(a)a$ |
| HJB-4 | $f(x, y, a) = F_2(x, y)m_{p_1(x,y)q_1(x,y)}(a)a$ |
| HJB-5 | $f(x, y, a) = F_4(x, y)m_{p_2(x,y)q_2(x,y)}(a)a$ |

, where $m_{\lambda, \mu}(a) = \left(1 + (\lambda a_1 + \mu a_2)^2\right)^{-\frac{1}{2}}$ for $\lambda, \mu \in R$ and for $c_1, c_2, c_3, c_4 > 0$,

$$C(x) = c_1 \sin\left(\dfrac{c_2 \pi x}{c_3} + c_4\right), c_1, c_2, c_3, c_4 > 0, \quad (F_1(x, y), F_2(x, y)) = \begin{cases} (0.5, 1) & y \geq C(x) \\ (2, 3) & y \leq C(x) \end{cases}$$

$$(F_3(x, y), F_4(x, y)) = \begin{cases} (0.2, 0.8) & y > C(x) + 0.25 \text{ or } y \leq C(x) - 0.25 \\ (1, 3) & C(x) < y \leq C(x) + 0.25 \\ (1, 1) & C(x) - 0.25 < y \leq C(x) \end{cases}$$

$$M_1(x, y) = \sqrt{\frac{\frac{F_2^2(x, y)}{F_1^2(x, y)}}{1+C'^2(x)}}, \quad q_1(x, y) = -M_1(x, y), \quad p_1(x, y) = M_1(x, y)C'(x)$$

$$M_2(x, y) = \sqrt{\frac{\frac{F_4^2(x, y)}{F_3^2(x, y)}}{1+C'^2(x)}}, \quad q_2(x, y) = -M_2(x, y), \quad p_2(x, y) = M_2(x, y)C'(x).$$

In all the following tests, we set the target $\Gamma = (0, 0)$ and $\Omega = [-0.5, 0.5] \times [-0.5, 0.5]$. We use the Computer with Core i 3, 2.53GHz for computing.

### 5.1 Convergence test

***Test 1.*** Here we solve equation HJB-1. In this case, the exact solution of HJB-1 is given by

$$u(x_1, x_2) = \sqrt{(1+\lambda^2)x_1^2 + (1+\mu^2)x_2^2 - 2\lambda\mu x_1 x_2} . \quad ([11])$$

We set $\lambda = 5, \mu = -10$. We denote the exact solution of HJB-1 as $u_0$.

We compare our numerical solution with the exact solution on various size of rectangular grids.

In Fig 5.1, we show the contour lines of $u_0$ and $u$ obtained using our method.

In table 5.1, we present the $L^1, L^2$ and $L^\infty$ norms of the error between the approximate solution obtained using our method and the exact solution. We also present the CPU times.

In table 5.1, $E_1 = \int_\Omega |u - u_0| \, dx, \ E_2 = \left( \int_\Omega |u - u_0|^2 \right)^{1/2}, \ E_\infty = \max_\Omega |u - u_0|$.

| $N_x = (N_y)$ | $E_\infty$ | $E_1$ | $E_2$ | CPU time |
|---|---|---|---|---|
| 101 | 0.026509 | 0.001085 | 0.002214 | 0.411 |
| 201 | 0.017757 | 0.000629 | 0.001380 | 1.757 |
| 401 | 0.012297 | 0.000420 | 0.000976 | 6.847 |

Table 5.1(HJB-1) The error between the exact solution and approximate solution obtained

our method and CPU time

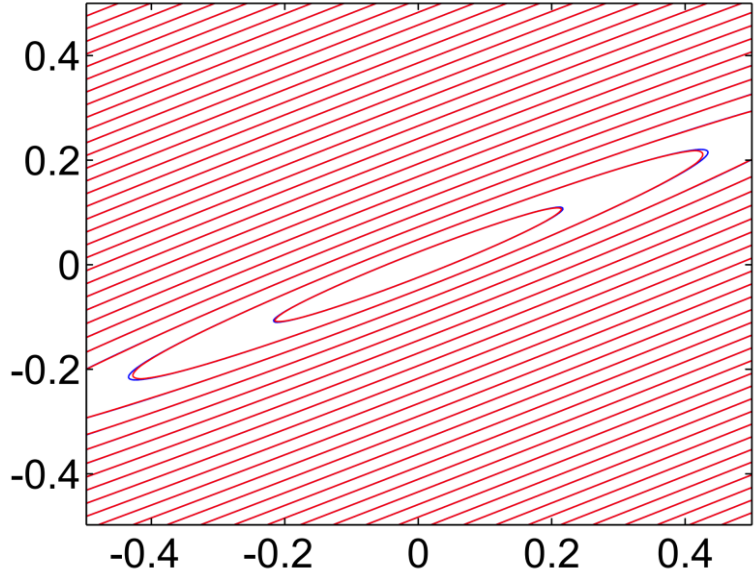

Fig. 5.1 (HJB-1, $N_x = N_y = 400$), Contour plots of $u_0$ (colored with red) and $u$ (colored with blue) obtained using our method

We have observed that our method captures viscosity solution of HJB-1 efficiently and accurately from fig. 5.1 and table 5.1.

## 5.2. Efficiency test

We also test our method with HJB-2~ HJB-4 to show the efficiency of our method.

For these equations, solution computed by OUM on a 801×801 grid will be referred to as the "exact" solution and will be denoted by $u^{exact}$.

***Test 2.*** Here we solve equation HJB-2 for $g(x, y) = 0.9\sin(2\pi x)\sin(2\pi y)$ as in [11]

On a manifold $z = g(x, y)$, the manifold's geodesic distance function $u$ is the viscosity solution of the above equation.

In Fig.5.2 and Table 5.2, we show the numerical result of our method for this equation.

| $N_x = (N_y)$ | $E_\infty$ | $E_1$ | $E_2$ | CPU time(s) |
|---|---|---|---|---|
| 100 | 0.074309 | 0.013985 | 0.021332 | 0.575 |
| 200 | 0.016033 | 0.003654 | 0.004721 | 2.200 |
| 400 | 0.003349 | 0.001302 | 0.001571 | 8.852 |

Table 5.2 (HJB-2) The error between the "exact solution" $u^{exact}$ and approximate solution

obtained our method and CPU time

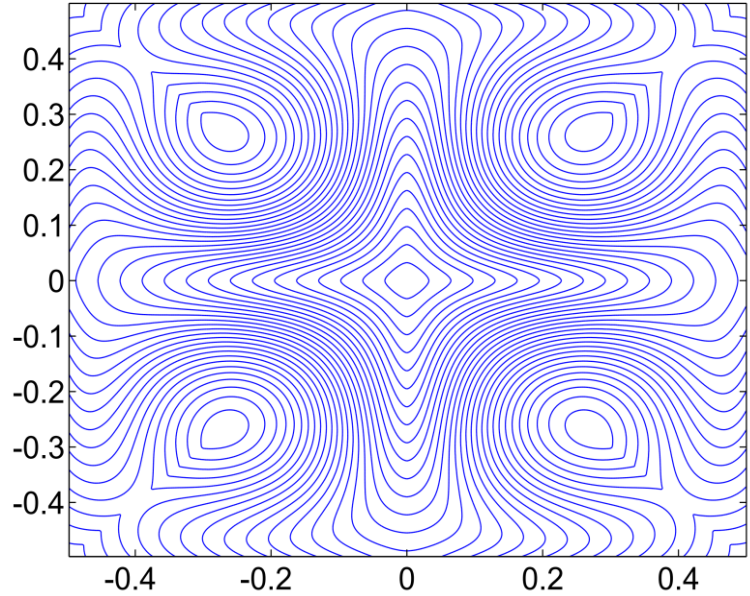

Fig. 5.2 (HJB-2, $N_x = N_y = 400$), Contour plots of $u$ (colored with blue) obtained using our method.

**Test 3.** Here we solve equation HJB-3 for $\lambda = 5, \mu = -10$ as in [1, 10].

In Fig.5.3 and Table 5.3, we show the numerical result of our method for this equation.

| $N_x = (N_y)$ | $E_\infty$ | $E_1$ | $E_2$ | CPU time(s) |
|---|---|---|---|---|
| 100 | 0.018783 | 0.007593 | 0.008857 | 0.427 |
| 200 | 0.012211 | 0.002445 | 0.003212 | 1.774 |
| 400 | 0.009819 | 0.002275 | 0.002843 | 6.946 |

Table 5.3 (HJB-3) The error between the "exact" solution $u^{exact}$ and approximate solution obtained our method and CPU time

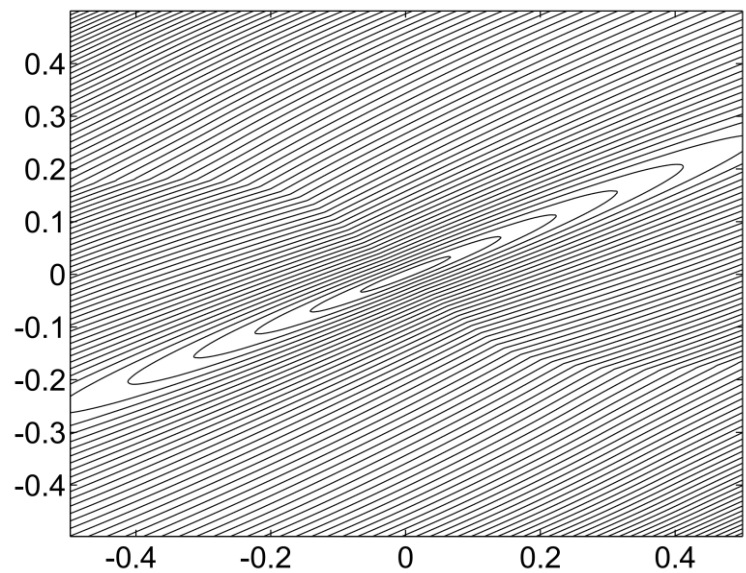

Fig. 5.3 (HJB-2, $N_x = N_y = 400$), Contour plots of $u$ (colored with black) obtained using our method.

**Test 4.** Here we solve equation HJB-4 for $\lambda = 5, \mu = -10$ as in [1, 10]

In Fig.5.4 and Table 5.4, we show the numerical result of our method for this equation.

| $N_x = (N_y)$ | $E_\infty$ | $E_1$ | $E_2$ | CPU time(s) |
|---|---|---|---|---|
| 100 | 0.023268 | 0.004361 | 0.006848 | 0.509 |
| 200 | 0.011557 | 0.001968 | 0.003044 | 2.053 |
| 400 | 0.004756 | 0.000614 | 0.000988 | 8.309 |

Table 5.4 (HJB-4) The error between the "exact" solution $u^{exact}$ and approximate solution obtained our method and CPU time

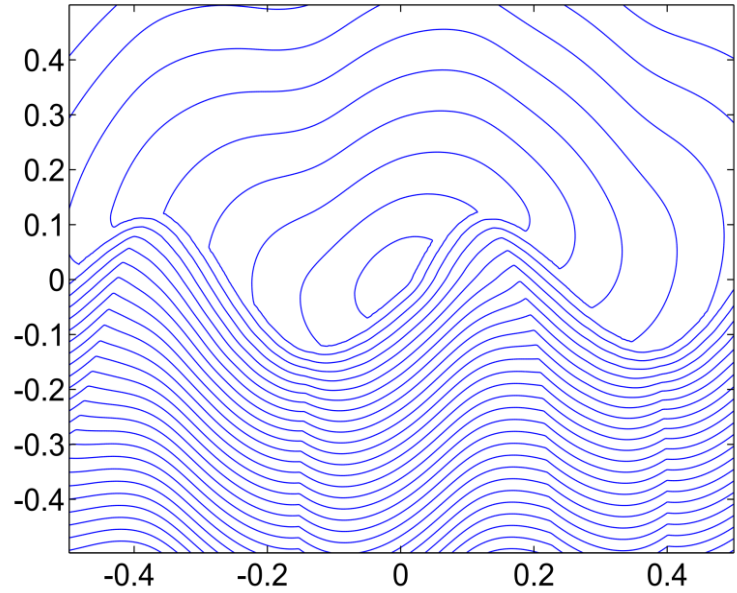

Fig. 5.4 (HJB-4, $N_x = N_y = 400$), Contour plots of $u$ (colored with blue) obtained using our method

***Test 5.*** Here we solve equation HJB-5 for $\lambda = 5, \mu = -10$ as in [11].

In Fig.5.5 and Table 5.5, we show the numerical result of our method for this equation.

| $N_x = (N_y)$ | $E_\infty$ | $E_1$ | $E_2$ | CPU time(s) |
|---|---|---|---|---|
| 100 | 0.05100 | 0.01815 | 0.016422 | 0.525 |
| 200 | 0.025654 | 0.005517 | 0.007107 | 2.085 |
| 400 | 0.012993 | 0.002192 | 0.002826 | 8.505 |

Table 5.5 (HJB-5) The error between the "exact" solution $u^{exact}$ and approximate solution obtained our method and CPU time

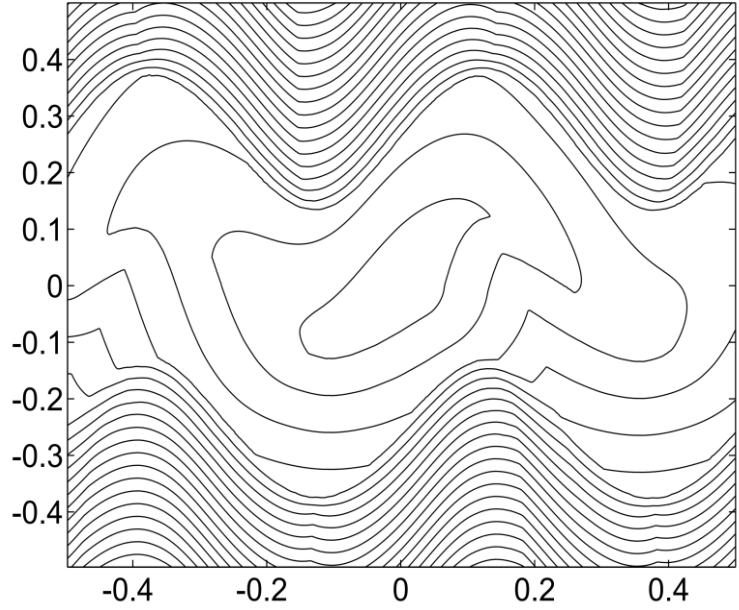

Fig. 5.5 (HJB-5, $N_x = N_y = 400$), Contour plots of $u$ (colored with black) obtained using our method.

## 6. Conclusion

We have proposed a new single-pass method for approximating the solution to an anisotropic eikonal equation related to the anisotropic min-time optimal trajectory problem. We have tested the new method with several anisotropic eikonal equations and have observed that it works well while significantly reducing the computational cost. The method presented in this paper is applicable for the anisotropic eikonal equations (1.3) and (1.4). We are currently studying on extending our technique to a wider class of problems including the generalized eikonal equation in a moving medium. We believe that our technique will be applicable for more general Hamilton-Jacobi-Bellman equations, which are associated to min-time problems.

## Appendix

**Lemma 1** *The viscosity solution $u(x)$ of (2.1) is semi-concave.*

*Proof.* Let us prove that $H(x, p)$ is Lipchitz continuous in $x$.

For arbitrary $x, y$, suppose that

$$H(x, p) = -\min_{a \in S_1}\{p \cdot a f(x, a)\} - 1 = -p \cdot a_0 f(x, a_0) - 1,$$

$$H(y, p) = -\min_{a \in S_1}\{p \cdot a f(y, a)\} - 1 = -p \cdot a_1 f(y, a_1) - 1.$$

If $p \cdot a_0 f(x, a_0) \geq p \cdot a_1 f(y, a_1)$, then

$$\begin{aligned}|H(x, p) - H(y, p)| &= |p \cdot a_0 f(x, a_0) - p \cdot a_1 f(y, a_1)| \\ &\leq |p \cdot a_0 f(y, a_0) - p \cdot a_0 f(x, a_0)| \leq \|p\| L \|x - y\|\end{aligned}.$$

If $p \cdot a_0 f(x, a_0) < p \cdot a_1 f(y, a_1)$, then

$$|H(x, p) - H(y, p)| = |p \cdot a_0 f(x, a_0) - p \cdot a_1 f(y, a_1)|$$
$$\leq |p \cdot a_1 f(x, a_1) - p \cdot a_1 f(y, a_1)| \leq \|p\| L \|x - y\|.$$

Therefore $H(x, p)$ is Lipchitz continuous in $x$. We also know that $H(x, p)$ is convex in the second argument and thus $u$ is semi-concave in $\Omega$.([31]) □

**Lemma 2** *Assume that $u$ is differentiable at $x \in \Omega$ and suppose that $\alpha(\cdot)$ is an optimal control for $x$ and $y(\cdot)$ is the optimal trajectory for $x$. Then $\nabla u(y(\cdot))$ is continuous at $[0, +\infty)$.*

*Proof.* For arbitrary $t > 0$, $u$ is differentiable at $y(t)$ because characteristics never emanate from the shocks-non differentiable point. ([11])

If $\nabla u$ is defined at $x \in \Omega \setminus \partial\Omega$, there exists constant $L$ such that $\|\nabla u\| \leq L$.([11])

We fix any $t \geq 0$. For arbitrary sequence $t_n \to t (t_n \geq 0)$, $y(t_n) \to y(t)$. The $\nabla u(y(t_n))$ is bounded. Therefore, it has convergent subsequence. Since $u$ is semi-concave, limits of all convergent subsequences are the same .i.e., $\nabla u(y(t))$. ([32]) Therefore $\nabla u(y(t_n)) \to \nabla u(y(t))$ and thus $\nabla u(y(\cdot))$ is continuous at $[0, +\infty)$. □

**Lemma 3** *Consider a grid node $x \in X$ and assume that $u$ is differentiable at $x$. Suppose that the characteristic for $x$ intersects the line segment $x_1 x_2$ at $x_0$, where $x_1, x_2 \in ND(x), x_1 \in N(x_2)$. Then*

$$\|\nabla u(x_0) - \nabla u(x)\| = O(h).$$

*Proof.* Suppose that $\alpha(\cdot)$ is an optimal control for $x$ and $y(\cdot)$ is the optimal trajectory for $x$.

If $y(t_0) = x_0$, since $\nabla u(y(\cdot))$ is continuous at 0, we get

$$\|\nabla u(y(t_0)) - \nabla u(y(0))\| = O(t_0).$$

By the Bellman's optimality principle, we get

$$t_0 \leq \frac{x_0 - x_1}{F_2}$$

and

$$\frac{x_0 - x_1}{F_1} \leq t_0.$$

Therefore $O(t_0) = O(h)$. From the fact that $y(0) = x, y(t_0) = x_0$, we conclude that

$$\|\nabla u(x_0) - \nabla u(x)\| = O(h) .□$$

**Lemma 4** *If $xx_1x_2$ is a sufficiently small simplex, which contains the characteristic for $x$, then*

$$\min(u(x_1), u(x_2)) < u(x).$$

*Proof.* By the Bellman's optimality principle, we get

$$u(x) = u(x_s) + \frac{\|x_s - x\|}{f(x, (x_s - x)/\|x_s - x\|)} + O(h^2) =$$
$$\lambda u(x_1) + (1-\lambda)u(x_2) + \frac{\|x_s - x\|}{f(x, (x_s - x)/\|x_s - x\|)} + O(h^2).$$

Since $xx_1x_2$ is sufficiently small simplex, we obtain that $u(x) > \lambda u(x_1) + (1-\lambda)u(x_2)$.

Therefore, $\min(u(x_1), u(x_2)) < u(x)$. □